\documentclass{amsart}
\usepackage{amsmath,amsfonts,amsthm,amscd,amssymb,color,comment,url}
\usepackage[mathscr]{eucal}
\usepackage[small,nohug,heads=LaTeX]{diagrams}
\usepackage[all,color]{xy}
\diagramstyle[labelstyle=\scriptstyle]

\newtheorem{theorem}{Theorem}[section]
\newtheorem{lemma}[theorem]{Lemma}
\newtheorem{proposition}[theorem]{Proposition}
\newtheorem{corollary}[theorem]{Corollary}

\newtheorem*{carlssonconj}{Conjecture (Carlsson)}

\newtheorem*{mainthm}{Theorem~\ref{maintheorem}}
\newtheorem*{propK_0xK*P}{Proposition~\ref{endomorphisms of irreds of k<G> is L}}
\newtheorem*{prophomotopygps}{Proposition~\ref{prop homotopy groups}}
\newtheorem*{propK0}{Proposition~\ref{K0 Repk<G> = Z[Z/linfZ]}}

\theoremstyle{remark}
\newtheorem{remark}[theorem]{Remark}

\theoremstyle{definition}
\newtheorem{definition}[theorem]{Definition}
\newtheorem{example}[theorem]{Example}

\DeclareMathOperator{\Br}{Br}

\DeclareMathOperator{\Gal}{Gal}
\DeclareMathOperator{\Hom}{Hom}
\DeclareMathOperator{\M}{M}
\DeclareMathOperator{\Mor}{Mor}
\DeclareMathOperator{\Stab}{Stab}

\DeclareMathOperator{\COLIM}{colim}
\DeclareMathOperator{\LIM}{lim}

\newcommand{\Tot}{\mathrm{Tot}}

\newcommand{\LRep}[2]{\mathrm{Rep}_{#1}[#2]}
\newcommand{\Rep}[2]{\mathrm{Rep}_{#1}\hspace{-.03in}\left\langle #2 \right\rangle}

\newcommand{\FF}{\mathbb{F}}
\newcommand{\NN}{\mathbb{N}}
\newcommand{\QQ}{\mathbb{Q}}
\newcommand{\ZZ}{\mathbb{Z}}

\newcommand{\XX}{\mathbb{X}}
\newcommand{\Z}{\mathscr{Z}}
\newcommand{\F}{\mathscr{F}}
\newcommand{\E}{\mathscr{E}}
\newcommand{\C}{\mathscr{C}}

\newcommand{\V}{\mathscr{V}}
\newcommand{\scr}[1]{\mathscr{#1}}
\renewcommand{\L}{\mathscr{L}}

\newcommand{\I}{\mathscr{I}}
\newcommand{\id}{\mathrm{id}}
\newcommand{\st}{\text{ s.t. }}
\newcommand{\xto}[1]{\xrightarrow{\,#1\,}}
\newcommand{\tensorover}[1]{{\otimes}_{#1}}

\newcommand{\newtensorover}[1]{{\otimes}_{#1}}

\newcommand{\smashover}[1]{\,{\wedge}_{#1}\,}

\newcommand{\leftexp}[2]{\,{\vphantom{#2}}^{#1\hspace{-.016in}}{#2}}

\newcommand{\EG}{E\langle G\rangle}

     \newcommand{\hatl}{^{^\wedge}_{\l}}

     \newcommand{\Hl}{\mathbb{H}_\l}
\newcommand{\T}{\mathscr{T}}
\newcommand{\Ti}[1][i]{\T_{#1}}
\newcommand{\Zmod}[1]{\ZZ/#1\ZZ}

     \newcommand{\Zl}{\ZZ_\l}

\newcommand{\G}{\mathscr{G}}
\newcommand{\N}{\mathscr{N}}

     \newcommand{\Gl}{\ZZ_\l\rtimes\ZZ_\l}
\newcommand{\dc}[2]{\left(#1\right)^{^\wedge}_{#2}}
\newcommand{\sdc}[2]{#1^{^\wedge}_{#2}}

     \newcommand{\sdcl}[1]{#1^{^\wedge}_{\alpha_\l}}

     \newcommand{\dcl}[1]{\left(#1\right)^{^\wedge}_{{\alpha_\l}}}

     \newcommand{\KTl}{\dcl{K\T}}
\newcommand{\KT}{{K\T}}
\newcommand{\RkG}{\Rep{k}{\G}}
\newcommand{\R}{\scr{R}}
\newcommand{\KR}{K\scr{R}}
\newcommand{\Gi}[1][i]{{\G}_{#1}}
\newcommand{\Ni}[1][]{\mathscr{N}_{i#1}}
\newcommand{\ki}[1][i]{k_{#1}}
\newcommand{\Di}[1][]{\Delta_{i#1}}
\renewcommand{\O}[1][]{\mathscr{O}_{#1}}
\newcommand{\Oi}[1][i]{\O[#1]}

\newcommand{\ROi}[1][i]{\Rep{\O[#1]}{\Gi[#1]}}
\newcommand{\Fi}[1][i]{\F_{#1}}
\newcommand{\RkGi}[1][i]{\Rep{\ki[#1]}{\Gi[#1]}}
\newcommand{\Ri}[1][]{\R_{i#1}}
\newcommand{\KRi}[1][]{K\Ri[#1]}
\newcommand{\KTi}{K\Ti}
\newcommand{\KZi}{K\Z_i}
\newcommand{\Zi}[1][]{\Z_i(#1)}

\newcommand{\colim}{\mathop{\COLIM}\limits}

\renewcommand{\lim}{\mathop{\LIM}\limits}

\renewcommand{\l}{\ell}
\newcommand{\EXACT}{\mathrm{ExactCategs}}
\newcommand{\ON}[1]{\mathrm{OpNm}\left({#1}\right)}
\newcommand{\NoTor}[1]{\mathrm{NoTor}\left({#1}\right)}
\newcommand{\Mod}[1]{\mathrm{Mod}\left({#1}\right)}

\newcommand{\tensorup}[2]{\left({#1} \tensorover{#2} -\right)}
\newcommand{\bracket}[1]{\langle #1 \rangle}

\newcommand{\op}{^{\mathrm{op}}}
\newcommand{\fform}{\mu}
\newcommand{\overbar}[1]{{\overline{#1}}}

\newcommand{\bigdca}[1]{\big({#1}\big)^{^\wedge}_{\alpha_\ell}}
\newcommand{\Fbar}{\overline{F}}
\DeclareMathOperator{\Vect}{Vect}
\newcommand{\smashup}[2]{\left({#1} \smashover{#2} -\right)}

\usepackage{xcolor}

\title{The $\ell$-adic $K$-theory of a $p$-local field}
\author{Grace K. Lyo (Massachusetts Institute of Technology)}

\begin{document}

\begin{abstract}
We verify a special case of a conjecture of G. Carlsson that describes the $\l$-adic $K$-theory of a field $F$ of characteristic prime to $\l$ in terms of the representation
theory of the absolute Galois group $G_F$. This conjecture is known to hold in two cases; in this article we examine the second case, in which the field in question is the field of Laurent series over a finite field $\FF_{p}((x))$.\end{abstract}
\maketitle

\section{Introduction}
\label{introduction}

We verify a special case of a conjecture of G. Carlsson that
describes the $\l$-adic $K$-theory of a field $F$ of
characteristic prime to $\l$ in terms of the representation
theory of the absolute Galois group $G_F$. This conjecture
employs a completion construction of Carlsson's called the {\em
  derived completion}~\cite{CarlDC}, which generalizes Bousfield
and Kan's $\l$-adic completion~\cite{Bk} to the categories of
ring and module spectra.  Throughout this article, we will work in
the category of $\mathbb{S}$-algebras (which we will refer to as simply ``ring
spectra'') from \cite{EKMM}. If $A$ is a commutative ring spectrum, $M$ is an $A$-module, and $\gamma\colon A\to B$ is a map of commutative ring spectra, then we may take the derived completion of $M$ with respect to $\gamma$; we will write $\sdc{M}{\gamma}$ or $\sdc{M}{B}$ to denote this completion.

We will be taking derived completions of $K$-theory spectra of
categories of representations. For any field $E$ equipped with
the action of a  group $G$, denote by $\LRep{E}{G}$ the category of finite dimensional continuous $E$-linear representations of $G$ and by $\Rep{E}{G}$ the analogous category of $E$-{\em semi}linear representations.
We define an $E$-semilinear representation of $G$ to be an
$E$-vector space $V$ with a $G$-action satisfying the relation
\[g(a\cdot v)=\leftexp{g}{a}\cdot g(v),\]
where $g\in G,\ a\in E,\ v\in V$, $\leftexp{g}{a}$ is the image of $a$ under the action by $g$ on $E$,
and $g(v)$ is the image of $v$ under the action by $g$ on $V$. It
was shown in \cite{EM} that the spectra $K \Rep{E}{G}$ are commutative
$\mathbb{S}$-algebras.

Denote by $\mathfrak{F}_\l$  the symmetric monoidal category with $\l$ objects,  only identity morphisms whose monoid structure is determined by the additive structure of $\FF_\l$. All derived completions will be taken with respect to a map
\begin{equation}\label{aug}
\alpha_\l\colon K \Rep{E}{G} \rTo K\mathbb{F}_\l\cong H\mathfrak{F}_\l ,
\end{equation}
for a prime $\l$ different from the characteristic of $E$.  This map is induced by the augmentation functor from $\Rep{E}{G}$ to $\mathfrak{F}_\l$, which sends a vector space to its $E$-dimension modulo~$\l$. 

Let $F$ be a field and $k$ an algebraically closed subfield of
$\overline{F}$.  Since $k$ is algebraically closed, it is
preserved by the absolute Galois group $G_F$ and hence there is a
category Rep$_k \langle G_F \rangle$.  Define
$\E_k^{\overline{F}}\colon$ $\Rep{k}{G_F}\rTo \Rep{\Fbar}{G_F}$
to be the extension of scalars functor, where the action by $G_F$
on $k$ and on $\overline{F}$ is the Galois action and where the
action by $G_F$ on $\E_k^{\overline{F}}(V)$ is the obvious one.
Denote by $KF{\hatl}$ the Bousfield--Kan $\ell$-adic completion.
\begin{carlssonconj}
For all primes $\l$ different from the characteristic of $F$, the induced map on completions (which we also denote $\E_k^{\Fbar}$)
\[
        \E_k^{\overline{F}}\colon \dc{K\Rep{k}{G_F}}{\alpha_\l} \rTo \dc{K\Rep{\overline{F}}{G_F}}{\alpha_\l} \simeq (KF)\hatl
\]
is a weak equivalence.
\end{carlssonconj}

\noindent Note that the weak equivalence on the right is  a
result of two facts.  First, the category $\Rep{E}{G}$ is
equivalent to the category $\Vect(E^G)$ of finite dimensional
$E^G$-vector spaces that arises whenever $G$ acts faithfully on
$E$.  Second, the derived completion $(KF)^{\wedge}_{\alpha_\ell}$
becomes the Bousfield--Kan $\ell$-adic completion when $\ell$ is
prime to the characteristic of $F$ \cite{CarlDRT}.

In~\cite{CarlDC}, Carlsson proved that this conjecture holds when $F$ is a characteristic zero field containing an algebraically closed subfield $k$, provided $G_F$ is abelian.

Fix primes $\l \ne p$, and for the remainder of this article, let
$k:=\overline{\FF}_p$ and define the field $\F$ to be the
compositum of all $p$ and $\l$-prime extensions of the field of
Laurent series $\FF_p((x))$.   Denote by $\F^t$ the maximal tame
extension of $\F$ and $G^t_{\F}$ the tame Galois group.  Our main goal is to prove the following theorem.

\begin{mainthm}
The map
\begin{equation}\label{main equivalence}
\E_k^{\F^t}\colon \dc{K\Rep{k}{G_\F^t}}{\alpha_\l} \rTo^{\simeq} \dc{K\Rep{\F^t}{G_\F^t}}{\alpha_\l}\simeq \dc{K\F}{\l},
\end{equation}
where $\F^t$ is the maximal tame extension of $\F$,
is a weak equivalence.
\end{mainthm}
\noindent For simplicity, we will write $\E=\E_k^{\F^t}$.

Note that in this theorem we use the tame Galois group
$G_\F^{t}$, which we denote by $\G$ throughout this article, and
the maximal tame extension $\F^t$, whereas the conjecture regards the absolute Galois group and the separable closure.

This substitution is justified by the fact that if $P$ is a
normal $p$-power subgroup of a profinite group $G$, and $G$ acts on a characteristic $p$ field $E$, then there is a weak equivalence
\[ K\Rep{E}{G/P}\simeq K\Rep{E}{G}.\]

Combining this with~\eqref{main equivalence} gives the conjectured equivalence $\E_{k}^{\overline{\F}}$.

The first step in our strategy towards proving this theorem will
be  to compute the $K$-groups $K_*\Rep{k}{\G}$,
(Section~\ref{homotopy groups}).  We plug this data in to the derived completion's algebraic to geometric spectral sequence, which converges to $\pi_*\big(\dc{K\Rep{k}{\G}}{\alpha_\l}\big)$.
In order to do this, we use results from ``Semisimple skew group
rings and their modules''~\cite{lyorib}, which is joint work with
K. Ribet. We  summarize the relevant portions of this article in
Section~\ref{skew group rings}. Applying these results yields the
following propositions about the homotopy groups of
$K\Rep{k}{\G}$.

\begin{propK_0xK*P}
The endomorphism rings of the irreducible representations in $\Rep{k}{\G}$ are matrix rings over the fixed field $k^\G$ and thus
\[
        K_*\Rep{k}{\G}\cong (K_0\Rep{k}{\G})\otimes K_*(k^\G).
\]
\end{propK_0xK*P}
\noindent In the sequel, define $\L=k^\G$.

\begin{propK0}
There is an isomorphism of rings, $K_0\Rep{k}{\G}\cong \ZZ[\QQ_\l/\ZZ_\l]$.
\end{propK0}

Finally, we compute the homotopy groups of $\dc{K\Rep{k}{\G}}{\l}$ and find that they are isomorphic to those of $\dc{K\F}{\l}$.

\begin{prophomotopygps}
The homotopy groups of $\dc{\Rep{k}{\G}}{\alpha_\l}$ are
\begin{align*}
\pi_n\big(\dc{K\Rep{k}{\G}}{\alpha_\l}\big)& \cong \pi_n\big(\dc{K\FF_p}{\l}\big)\oplus \pi_n\Sigma\big(\dc{K\FF_p}{\l}\big)\\
&\cong \pi_n\big(\dc{K\F}{\l}\big)\, .
\end{align*}
\end{prophomotopygps}

In Section~\ref{weak equiv}, we show that the map $\E$ induces the isomorphism of homotopy groups by factoring $\E$ through the intermediate category
\[
\colim_{\G_i' \in \ON \G} \Rep{\Oi}{\G/\G_i'},
\]
where $\ON \G$ is defined to be the opposite category of the
category of open normal subgroups of $\G$, $\Oi$ is the ring of
integers in the fixed field $\F^{\G_i'}$, and
$\Rep{\Oi}{\G/\G'_i}$ is the category of finitely generated
$\Oi$-modules equipped with the continuous semilinear action by
$\G/\G'_i$ on $\F^{\G'_i}$.
\\

{\bf Acknowledgments.} I would like to thank Tyler Lawson for his many helpful suggestions on a draft of this paper.  

\section{Completions}
\label{background}

\subsection{The Bousfield-Kan Completion}
We briefly review the Bousfield-Kan completion~\cite{Bk}. Let $\C$ be a category.

\begin{definition}
A monad (or triple) on $\C$ is a functor $T \colon \C \to \C$ together with natural transformations $\eta \colon \id \to T$ and $\mu \colon T \circ T \to T$ such that for any object $X \in \C$ the diagrams below commute.
\begin{center}
\begin{tabular}{cc}
\begin{diagram}
TX&\rTo^{T\eta(X)}&T \circ TX & \lTo^{\eta(TX)} & TX\\
&\rdTo_{=}&\dTo_{\mu(X)}&\ldTo_{=}\\
&&TX
\end{diagram} &
\begin{diagram}
T^3(X)&\rTo^{\mu(TX)}&T^2(X)\\
\dTo^{T\mu(X)}&&\dTo_{\mu(X)}\\
T^2(X)&\rTo^{\mu(X)}&T(X)
\end{diagram}
\end{tabular}
\end{center}
\end{definition}
See~\cite{May} for more details.

\begin{definition}
Suppose that $c\in \C$. The cosimplicial resolution of $c$ relative to $T$ is the cosimplicial object $\R^\bullet(T,c)$ whose $n$ simplices are obtained by making $n+1$ applications of $T$. That is, for any nonzero integer $n$,
\[\R^n(T,c):=T^{n+1}(c).\]
The coface and codegeneracies are determined by $\eta$ and $\mu$ respectively:
\[\delta^s \colon T^{k+1}(c) \to T^{k+2}(c) = T^s (\eta(T^{k+1-s}(c))\]
and
\[\sigma^t \colon T^{k+1}(c) \to T^k (c)=T^t(\mu(T^{k-t-1}(c)).\]
\end{definition}

Now let $\C$ be the category of simplicial sets and $R$ a
ring. Define $F_R$ to be the composition of the forgetful functor
with the free $R$-module functor applied levelwise. Then $F_R$ (together with the obvious natural transformations) is a monad on $\C$.

\begin{definition}[Bousfield-Kan]
Let $X_\bullet$ be a simplicial set. Then the $R$-completion of $X_\bullet$, denoted $R_\infty X$, is defined to be the cosimplicial resolution of $X_\bullet$ relative to $F_R$. The $\l$-adic completion of $X_\bullet$, denoted $X\hatl$, is the $R$-completion for $R=\FF_\l$.
\end{definition}

This construction can be extended to the category of spectra.

\begin{definition}
Suppose that $X=\{X_i\}$ is a spectrum. Then the $R$-completion of $X$, denoted $R_\infty X$ is defined to be the ``spectrification'' of the prespectrum $\{R_\infty(X_i)\}$. As before, the $\FF_\l$-completion is called the $\l$-adic completion and is denoted $X\hatl$.
\end{definition}

The following propositions are about the homotopy groups of
$\l$-adically completed spectra.

\begin{proposition}\label{pi of completion}
Let $X$ be a spectrum whose homotopy groups are all finitely generated. Then
\[\pi_n(X\hatl)\cong\ZZ_\l\otimes\pi_nX.\]
\end{proposition}

\begin{proof}
This follows from Example~5.2 in \cite{Bk}.
\end{proof}

\begin{proposition}\label{l-comp(KL)=l-comp(KFl)}
Let $F$ be a finite field of characteristic $p$  and $L\supseteq F$ an extension of degree prime to $\l$. The inclusion $F\hookrightarrow L$ induces a weak equivalence,
\[(KF)\hatl \rTo ^{\cong} (KL)\hatl.\]
\end{proposition}

\begin{proof}
Let $L$ be a finite extension of $F$ of degree prime to $\l$. The composition of extension and restriction of scalars,
\begin{diagram}
KF &\rTo& KL &\rTo& KF
\end{diagram}
induces multiplication by $[L:F]$ on $\pi_*(KF)$. Since $[L:F]$ is prime to $\l$, the $\l$-primary components of $\pi_* KF$ are direct summands of the homotopy groups $\pi_*KL.$

Because $\pi_n KF$ and $\pi_n KL$ are both cyclic, they must be isomorphic.

\end{proof}

\begin{proposition}
The $\l$-adic $K$-groups of $\F$ are
\[
        \pi_*\dc{K\F}{\l}\cong \pi_*\big(\dc{K\FF_p}{\l}\big)\oplus \pi_*\Sigma\big(\dc{K\FF_p}{\l}\big).
\]
\end{proposition}

\begin{proof}

This is proven in Section~4.2 of~\cite{CarlDRT}.
\end{proof}

\subsection{The Derived Completion}

To define the derived completion, we observe that when $\gamma \colon A\rightarrow B$ is a map of commutative ring spectra, then the functor $\smashup{B}{A}$ together with the natural transformations
\[  \eta \colon \id \rightarrow \smashup{B}{A}, \]
defined by $\eta=f\wedge \id$ and
\[ \mu\colon (B\wedge_A B\wedge_A -) \rightarrow \smashup{B}{A}, \]
defined by multiplication on $B$, is a monad on the category
$\Mod{A}$ of $A$-modules. When $A$ is the sphere spectrum and $B$
is the Eilenberg-MacLane spectrum of a ring $R$, the resulting
monad $\big(H(R)\smashover {\mathbb{S}}-\big)$ is equivalent in
the category of spectra to the monad $F_R$. Thus for an
$\mathbb{S}$-module $M$ there is a weak equivalence of
$\mathbb{S}$-modules,
\[ \sdc{M}{HR}\simeq R_\infty M.\]
Denote by $\tilde{B}$ the functorial fibrant replacement for $B$ in the category of commutative $A$-algebras, and by $\hat{\mathscr{R}}^\bullet \big(\,(\tilde{B}\smashover A-),M\big)$ the functorial fibrant replacement for $\mathscr{R}^\bullet\big(\,(\tilde{B}\smashover A-),M\big)$ in the category of cosimplicial $A$-modules.

\begin{definition}
The derived completion $\sdc{M}{\gamma}$ of $M$ with respect to the map $\gamma \colon A\rightarrow B$ is the total space,
\[ \sdc{M}{\gamma} := \Tot \big(\hat{\scr{R}}^\bullet (\tilde{B}\smashover A-,\,M)\big). \]
\end{definition}

The following propositions and theorems from~\cite{CarlDC} are the most relevant for our purposes. We refer the reader to~\cite{CarlDC} for more details.

\begin{proposition}
Let $A\rightarrow B$ be a map of commutative ring spectra and $M$ an $A$-module. Then the following statements all hold.
\begin{enumerate}
\item The construction $M\rightarrow\sdc{M}{B}$ is functorial for homomorphisms of $A$-modules. The map on $B$-completions induced by a homomorphism $f\colon M\rightarrow N$ of $A$-modules will be denoted by $\sdc{f}{B}$.
\item Let $f\colon M\rightarrow N$ be a homomorphism of $A$-modules, where $A$ is a commutative ring spectrum,
and suppose that $\id_B \smashover A f$ is a weak equivalence of $B$-module spectra. Then the map on completions $\sdc{f}{B}$ is also a weak equivalence of spectra.
\item Let $M''\rightarrow M\rightarrow M'$ be a cofibration sequence of $A$-module spectra.
Then the sequence $\dc{M''}{B}\rightarrow \dc{M}{B}\rightarrow \dc{M'}{B}$ is a cofibration sequence up to homotopy.
\item Let $A''\rightarrow A\rightarrow B$ be a diagram of commutative ring spectra, and let $M$ denote a $A$-module spectrum. Denote by $M_{A^{''}}$ the module $M$ regarded as an $A''$-module. Suppose that the natural map $B\smashover {A''} (M_{A^{''}})\rightarrow B\smashover{A} M$ is a weak equivalence of spectra. Then the natural map $\dc{M_{A^{''}}}{B}\rightarrow\sdc{M}{B}$ is also an equivalence of spectra.
\item There is a natural map $\eta\colon M\rightarrow \sdc{M}{B}$.
\item Suppose $M$ is an $A$-module spectrum for which the $A$-module structure admits a $B$-module structure extending the given $A$-module structure. Then the natural map $\eta\colon M\rightarrow \sdc{M}{B}$ is an equivalence of spectra.
\end{enumerate}
\end{proposition}

\begin{theorem}\label{depends on pi_0}
Suppose that we have a diagram $A\rightarrow B\rTo^{\beta} B'$ of
commutative ring spectra and a left $A$-module spectrum~$M$. Suppose further that $A$, $B$, and $M$ are all $(-1)$-connected, that the homomorphism $\pi_0(\beta)$ is an isomorphism, and that the natural homomorphisms $\pi_0(A)\rightarrow \pi_0(B)$ and $\pi_0(A)\rightarrow \pi_0(B')$ are surjections. Then  the natural homomorphism $\sdc{M}{B}\rightarrow \sdc{M}{B'}$ is an equivalence of $A$-module spectra.
\end{theorem}

\begin{theorem}[Algebraic to geometric spectral sequence]\label{alg to geom}
  Let $A\rightarrow B$ be a map of commutative ring spectra, with
  $\pi_0 A\rightarrow \pi_0 B$ surjective, and let $M$ be an
  $A$-module spectrum. The ring $\pi_0A$ is commutative, the ring
  $\pi_0B$ is a commutative $\pi_0A$-algebra, and for each $i$,
  the abelian group $\pi_iM$ has a $\pi_0A$-module structure. We
  may therefore construct the derived completion
  $\dc{\pi_iM}{\pi_0B}$ for each $i$. There is a second quadrant
  spectral sequence with
  $E^{s,t}_1=\pi_{s+2t}\big(\dc{\pi_{-t}M}{\pi_0B}\big)$,
  converging to $\pi_{t+s} (\sdc{M}{B})$.
\end{theorem}

The following example from~\cite{CarlDC} will play an important role in our computation of the homotopy groups $\pi_* \big(\dc{K\Rep{K}{\G}}{\alpha_\l}\big)$.

\begin{example}\label{dc of Z/linfZ}
Denote by $R[\ZZ_\l]$ the ring of isomorphism classes of complex representations of the group of $\l$-adic integers. Then $R[\ZZ_\l]$ can be expressed as a colimit
\[ R[\ZZ_\l]=\colim_n R[\ZZ/\l\ZZ]. \]
For each finite representation ring, there is an isomorphism,
\[ R[\ZZ/\l^n \ZZ] \cong \ZZ[\mu_{\l^n}], \]
and thus $R[\ZZ_\l]\cong \ZZ[\mu_{\l^\infty}]$, where $\mu_{\l^\infty}$ is the group of all $\l$-power roots of unity.

Define $\alpha_\l$ as in~\eqref{aug}. Then the homotopy groups of the derived completion of the Eilenberg-MacLane spectrum $(HR[\ZZ_\l])^\wedge_{\alpha_\l}$ relative to $\alpha_\l$ are

\[ \pi_{*} \big(\dc{HR[\ZZ_\l]}{\alpha_\l}\big) \cong
    \begin{cases}
        \ZZ_\l & \textrm{for } i=0\textrm{ or }1,\\
        0 & \textrm{otherwise.}
        \end{cases} \]
\end{example}

\section{Semisimple skew group rings}
\label{skew group rings}

 In this section we summarize the relevant results of ``Semisimple skew group rings and their modules''~\cite{lyorib}, which presents a means of describing the objects, endomorphisms, and tensor products in the category of finitely generated modules over a skew group ring. We will compute the ring $K_0$ of isomorphism classes of finitely generated modules over a skew group ring under the operations $\oplus$ and $\otimes_E$.

 Suppose that $E$ is a field equipped with the action of a finite group $G$. Then we define the skew group ring $E\bracket{G}$ to be the $E$-vector space with basis set $G$ endowed with the multiplication structure determined by {$(\alpha g) (\beta h)=\alpha\hspace{.013in}\leftexp{g\hspace{-.01in}}{\beta}\,gh$}, for $\alpha$ and $\beta\in E$, for $g$ and $h\in G$, and where {$\leftexp{g\hspace{-.01in}}{\beta}$} is the image of $\beta$ under~$g$. If $G$ acts trivially on $E$, then $E\bracket{G}$ is the ordinary group ring $E[G]$. In our motivating example, the action by $G$ on $E$ is not trivial; however it is typically not faithful either, and hence the elements of $G$ that act trivially on $E$ form a proper normal subgroup $N$.

 Throughout this section, we will assume that $N$ is abelian and
 therefore that the action by $G$ on $N$ factors through the
 quotient $G/N$, which we will denote by~$\Delta$. The assumption
 that $N$ is abelian also implies that there is a cohomology
 group $H^2(G/N,N)$.
 The group $G$ is an extension of $G/N$ by $N$ and thus
 determines an equivalence class in $H^2(G/N,N)$.

 We will also assume that the order of $G$ is invertible in $E$
 and thus, by a generalization of Maschke's Theorem, that
 $E\bracket{G}$ is semisimple. Every $E\bracket{G}$ module can
 therefore be broken down as a direct sum of simple modules.
 Every simple module  is isomorphic to a simple module over a single factor
 in the Artin-Wedderburn decomposition of $E\bracket{G}$ as a
 product of matrix rings 
\[
        E\bracket{G}=\prod_{i=1}^sM(d_i,D_i),
\]
where $s$ and $d_i$ are positive integers and  $D_i$ is a division ring (which will contain the fixed field $F:=E^G$). This decomposition is unique; therefore, describing the objects and endomorphisms of $E\bracket{G}$-modules is equivalent to finding the division rings $D_i$ and the corresponding dimensions $d_i$. Each $D_i$ is a central simple algebra over its center $Z_i$, so $D_i$ corresponds to an element of $\Br(Z_i)$.

The product $Z=\prod Z_i$ is the center of $\EG$.

\begin{proposition}
The map
\[
        \mu\colon E \otimes_F Z\rTo E[N]
\]
defined by
\[
        e\otimes z \mapsto ez,
\]
where $e \in E$ and $z \in Z$, is an $E$-algebra isomorphism.
\end{proposition}

Now define $\Br(Z):=\prod\Br(Z_i)$. By a theorem of Auslander and Goldman~\cite{AusGold} we can identify $\Br(Z)$ with the cohomology group $H^2(G_F,(\Fbar\otimes Z)^*)$. The skew group ring $\EG$ thus corresponds to a 2-cocycle in $H^2(G_F,(\Fbar\otimes Z)^*)$. There is a natural map obtained by composing the maps
\begin{align*}
H^2(G/N,N)&\rTo^{\cong}H^2(\Gal(E/F),N)\\
          &\rTo^{\rm inf.}H^2(G_F,N)\\
          &\rTo^{\rm inc.}H^2(G_F,F[N]^*)\\
          &\rTo^{\mu}H^2(G_F,(\Fbar\otimes Z)^*)\\
          &\rTo^{\cong}\Br(Z),
\end{align*}
where the map ``{inf.}'' is inflation and ``{inc.}'' is inclusion. Call this composite~$\Phi$.

\begin{theorem}\label{thm Phi sends G to EG}
The map $\Phi$ sends $[G]\in H^2(G/N,N)$ to $[\EG]\in\Br(Z)$.
\end{theorem}

\begin{corollary}\label{cor semidirectproduct gives matrix ring}
If $\Br(F)=1$, there is an isomorphism of rings,
\[\EG\cong \M([E:F],Z).\]
\end{corollary}

\begin{proposition}\label{thm EG free over Z}
As a $Z$-module, $\EG$ is free of rank $[E:F]^2$.
\end{proposition}
\noindent This describes $\EG$ up to isomorphism.

To compute the multiplicative structure of $K_0(\EG)$, we first observe that the set $\Hom_F(Z,\overline{F})$ of $F$-algebra homomorphisms is endowed with a natural action by the absolute Galois group $G_F$ which is defined by
\[ (g*\phi)(z):= g(\phi(z))\]
for $g\in G_F$ and $\phi\in\Hom_F(Z,\overbar{F})$. This action partitions $\Hom_F(Z,\overbar{F})$ into a disjoint union of orbits $\Hom_F(Z_i,\overbar{F})=\coprod \Hom(Z_i,\overbar{F})$.

\begin{corollary}\label{Hom(Z,F-bar)=Hom(N,F-bar*)}
The map $\fform$ induces a bijecton of sets 
\[ \mu_*\colon \Hom (Z,\overline{F}) \rTo^\cong \Hom(N,\overline{F}^{*}) \]
that is $G_F$ equivariant if the action by $G_F$ on $\Hom(N,\overline{F}^{*})$ is defined in the natural way,
\begin{equation}
    (g * \phi)(n):= g(\phi(\leftexp{\overline{g}^{-1}}{n})), \label{eqtn skew G_F action}
\end{equation}
where $g\in G_F$, $\phi\in\Hom(N,\overline{F}^{*})$, and $n\in N$; the element $\overline{g}$ is the image of $g$ in the quotient $G/N\cong \Gal(E/F)$ of $G_F$ and $\leftexp{\overline{g}^{-1}}{n}$ is the image of $n$ under conjugation by $\overline{g}^{-1}$.
\end{corollary}

Under $\mu_*$, the orbits $\Hom(Z_i,\overline{F})$ get sent to orbits of $\Hom(N,\overline{F}^{*})$ which we call $\mathfrak{O}_i$; each of the $\mathfrak{O}_i$ therefore corresponds to a factor $ M(d_i,D_i)$ of $\EG$.
Note that the field $Z_i$ is isomorphic to a fixed field of $\overline{F}$ under the subgroup $\Stab_{G_F}(\chi)$ of elements in $G_F$ that fix any character $\chi$ in $\mathfrak{O}_i$. The degree of $Z_i$ over $F$ is equal to the number of elements in $\mathfrak{O}_i$.

Now fix a common separable closure $\overline{F}$ of $F$ and all separable extensions of $F$, and write $N^\vee:=\Hom(N,\Fbar^*)$. We will associate to every $\EG$-module $W$ an element of the group ring
$\ZZ[N^{\vee}]$ as follows.

We start by restricting $W$ to an $E$-linear representation of $N$. Since $N$ is abelian, we can write $\Fbar\tensorover{E}W$ as a sum of $1$-dimensional representations
\[\overline{F}\newtensorover{E}W\cong\bigoplus_{N^{\vee}} U_{\chi}^{\ m_\chi},\]
where $m_{\chi}\in\ZZ^{\geq 0}.$
This decomposition is unique and suggests a map $\rho$ defined by
\begin{align}\label{eqtn defn of rho}
\rho(W):=\sum m_{\chi}\chi
\end{align}
from $K_0\EG$ to the group ring $\ZZ[N^{\vee}]$.

\begin{proposition}\label{prop irrep stable sum}
For every $\EG$-module $W$, the image $\rho(W)$ is stable under the action by $G_F$ defined in equation~\eqref{eqtn skew G_F action}.
\end{proposition}

\begin{proposition}\label{prop irreducibles corresp to one orbit}
For a simple module $V_i$, the set of elements that appear with a nonzero coefficient in $\rho(V_i)$ is exactly the orbit $\mathfrak{O}_i$ of $N^{\vee}$, defined in Corollary \ref{Hom(Z,F-bar)=Hom(N,F-bar*)}. The map $\rho$ is therefore injective.
\end{proposition}

If $V_i$ is a simple module over the $i$th factor $M(d_i,D_i)$ of $\EG$, then the center $Z_j$ of a different factor, where $1\leq j \leq s$ and $j\neq i$, acts trivially on $V_i$. The elements of $Z_i$, on the other hand, act through the characters in $\Hom(Z_i,\overline{F})$ by
\[ z\cdot v=\chi(z)v,\]
where $z\in Z_i$, $v\in V$, and $\chi \in\Hom(Z_i,\overline{F})$. By Corollary \ref{Hom(Z,F-bar)=Hom(N,F-bar*)}, these characters correspond exactly to those in $\mathfrak{O}_i$.

It follows immediately from Propositions \ref{prop irrep stable sum} and \ref{prop irreducibles corresp to one orbit} that for a simple module $V_i$, there is some positive integer $m_i$, which we call the ``multiplicity'' of $V_i$, such that
\[\rho(V_i)=m_i\sum_{\chi\in\mathfrak{O}_i}\chi.\]

\begin{corollary}
If $\Br(F)=1$, then $m_i=1$ for all irreducible $\EG$-modules $V_i$.
\end{corollary}
\begin{theorem}
The map $\rho$ is multiplicative with respect to $\otimes_E$, and hence an $\EG$-module $W$ is irreducible exactly when its image in $\ZZ[N^{\vee}]$ is
\[\rho(W)=m_i\sum_{\chi\in\mathfrak{O}_i}\chi,\]
for some $i$.
If $W'$ is another module over $\EG$, then
\[W\tensorover{E}W'\cong\bigoplus_{\mathscr{V}} V_i^{r_i}\]
where $r_i\in\ZZ^{\geq 0}$ and $\V$ is a set containing one
representative $V_i$ for each isomorphism class of simple
modules, exactly when $\rho(W)\cdot \rho (W')=\sum_{\mathscr{V}}\rho(V_i^{r_i})$.
\end{theorem}

\begin{example}\label{K0 is Z[N]}
Suppose that the Brauer group $\Br(F)$ is trivial and therefore that $\Hom(N, \overline{F}^*)$ breaks into singleton orbits under the action by $G_F$. Then the map $\rho$ induces an isomorphism of rings
\[\rho:K_0 (E\bracket{G}) \rTo^{\cong} \ZZ\big[N^\vee\big].\]
\end{example}

\begin{example}\label{L<Z/pi x Z/pi>}
  We will show in Proposition~\ref{prop homotopy
    groups} that $\G\cong \Zl \rtimes_\xi \Zl$
  where the action is the Frobenius action.  Let
  $\G_i=\Zmod{\l^i} \rtimes \Zmod{\l^{i'}}$ be the quotient of
  the group $\G$ by the subgroup $\G_i'=\l^i \ZZ_\l \rtimes
  \l^{i'} \ZZ_\l$, where $i>i'$ are positive integers. The
  integer $i'$ is defined to be the degree of the Galois
  extension $\L(\zeta)/\L$, where $\zeta$ is a primitive
  $\l^{i}\mathrm{th}$ root of unity. The group $\G_i$ acts on $(\overline{\FF}_p)^{\l^i \ZZ_\l}$
  (which we denote $k_i$) in the obvious way, and the fixed field
  of this action is $\L$.

Let $\N_i:=\Zmod{\l^i}$ be the subgroup of elements of $\G_i$ that act trivially on $k_i$ and define $\Delta_i$ to be the quotient $\Zmod{\l^{i'}}\cong\G_i/\N_i$. Let $\eta$ be a generator of $\N_i$. The quotient group $\Delta_i$ (which is also a subgroup of $\G$) is generated by the Frobenius automorphism $\phi$ and sends an automorphism $\eta^r \in \N_i$, where $r$ is an integer between 0 and $\l^i-1$, to its $p\mathrm{th}$ power; that is, $\phi \eta^r \phi^{-1}=\eta ^ {p r}$. The action by the absolute Galois group $G_\L$ on $\N_i^\vee=\Hom(\N_i, \overline{\L}^*)$ divides $\N_i^\vee$ into singleton orbits.
By Example~\ref{K0 is Z[N]},
\[K_0\Big(k_i\big<\Zmod{\l^i} \rtimes \Zmod{\l^{i'}}\big>\Big) \cong \ZZ[\mu_{\l^i}].\]
\end{example}

Now we consider the case where $G$ is a profinite group. Define $\ON{G}$ to be the category whose objects are the open normal subgroups of $G$. For every pair $G''$ and $G'$ of such subgroups, the morphism set $\Mor(G'',G')$ has exactly one element if $G'$ is a subgroup of $G''$ and is empty otherwise. Let $\ON{G}\op$ be the opposite category. The functor
\[\Rep{E}{G/-}: \ON{G}\op \to \EXACT\]
takes an open normal subgroup $G' \subseteq G$ to the exact category $\Rep{E'}{G/G'}$, where $E'=E^{G'}$.

\begin{proposition}\label{k<G> is a colimit}
The category $\Rep{E}{G}$ is a colimit
\[\Rep{E}{G} \cong \colim_{G' \in \ON{G}\op} \Rep{E'}{G/G'}.\]
Each constituent $\Rep{E'}{G/G'}$ of the colimit may be identified with a subcategory of $\Rep{E}{G}$ under the functor
\[\varphi': \Rep{E'}{G/G'} \to \Rep{E}{G},\]
which is defined by
\[\varphi'(W)=E \tensorover{E'} W.\]
The action by an element $g \in G$ on an element $e \otimes w \in E \otimes W$ is defined by
\[g(e \otimes w)=\leftexp{g}{e} \otimes \overline{g}(w),\]
where $\leftexp{g}{e}$ is the image of $e$ under the action by $g$, the element $\overline{g}$ is the image of $g$ under the quotient map $G \to G/G'$, and $\overline{g}(w)$ is the image of $w$ under the action by $\overline{g}$.
\end{proposition}

\begin{corollary}
There is an isomorphism of rings,
\[
K_0\Rep{E}{G} \cong \colim_{G'\in \ON{G}} K_0E'\bracket{G/G'}.
\]
\end{corollary}

\section{The homotopy groups of $\dc{K\Rep{k}{\G}}{\alpha_\l}$}
\label{homotopy groups}

In this section we compute the homotopy groups of $K\Rep{k}{\G}$,
which we will feed in to the algebraic to geometric spectral sequence to obtain the homotopy groups of the derived completion $\dc{K\Rep{k}{\G}}{\alpha_\l}$.

Define $\NN$ to be the filtered category whose objects are the natural numbers and whose morphism sets $\Mor_{\NN}(m,n)$ consist of one element when $n$ divides $m$ and are empty otherwise. For any ring $R$ and any positive integer $r$, we will denote by $R[[x^{(1/r)},x^{(1/s)}]]$ the colimit,
\[R[[x^{(1/r)},x^{(1/s)}]]:=\colim_{i\in \NN \st r,s \nmid i} R[[x^{1/i}]].\]
If $R$ is already a colimit of the form
\[R=\colim_{i \in \I} R'[[x^{1/i}]]\]
for some ring $R'$ and some subset $\I \subseteq \NN$, then we denote by $R[x^{(1/r)}]$ the colimit
\[R[x^{(1/r)}]:=\colim_{i \in \I \st r \nmid i} R'[[x^{1/i}]].\]
We define $F((x^{(1/r)}))$ and $F(x^{(1/r)})$ in the analogous way.

We start by describing the tame Galois group $\G$. It is easy to see that there is a normal subgroup of $\G$, which we denote by $\N$, that is isomorphic to $\Gal(\overline{\L}/\L)\cong \Zl$.
The quotient by this normal subgroup, which permutes the roots of $x$, will be denoted $\Delta$.

\begin{proposition}\label{prop homotopy groups} 
The tame Galois group of $\F$ is the semidirect product $\Gl$.
\end{proposition}
\begin{proof}
The quotient $\Delta$ sits inside $\G$ as the Galois group $\Delta \cong \Gal\big(\F^t/\F^t((x^{(1/\l)}))\big)$. The subgroups $N$ and $\Delta$ intersect trivially and their product fixes exactly $\F$. Therefore $\G$ is the product of subgroups $\N\cdot\Delta$.

\end{proof}

\begin{remark}
The generator of $\Delta$ (the Frobenius automorphism) acts on $\N$ by raising elements to the $\l$-th power.
\end{remark}

\begin{remark}
The subgroup $\Gl$ is not necessarily normal in $G_{\FF_p ((x))}^t$. This is because for all primes $q$ different from $\l$ and $p$, $\F$ contains all $q$th roots of $x$ but not necessarily all such roots of $1$. For example, if $q\mid p^\l-1,$ then the order of $p$ modulo $q$ will be $\l$ and hence the $q$th roots of unity can only be adjoined by making a degree $\l$ extension. The extension $\F/\FF_{p}((x))$ is therefore not necessarily Galois.\end{remark}

Define $\G_i,\ \G_i',\ \N_i$, and $k_i$ as in Example~\ref{L<Z/pi x Z/pi>}. By Proposition~\ref{k<G> is a colimit} there is an equivalence of categories
\[
        \Rep{k}{{\G}}\cong \colim_{\G_i' \in \ON{\G}} \Rep{k_i}{{\G/\G_i'}}.
\]
\begin{proposition}\label{endomorphisms of irreds of k<G> is L}
The ring of endomorphisms of every irreducible object in $\RkG$ is $\L$. The $K$-groups of $\RkG$ are therefore
\[K_*\RkG\cong \big(K_0\,\RkG\big)\tensorover{} K_*\L.\]
\end{proposition}
\begin{proof}
The orbits of $\Hom(\Ni, k^*)$ under the action by the absolute Galois group $G_{\L}$ are all singleton orbits, and thus by Proposition~\ref{prop irreducibles corresp to one orbit}, there is a one-to-one correspondence between characters of $\Ni$ and irreducible representations of $\Gi'$. An irreducible representation $V$ that corresponds to the character $\chi \in \Hom(\Ni, \overline{k_i}^*)$ is a $k_i$-vector space on which $\Ni$ acts by
\[ n\cdot v =\chi(n)v,\]
for $n \in \Ni$ and $v \in V$.
 Since $\Gi$ is a semidirect product of $\Ni$ and $\Di$, every $g\in \Gi'$ can be written uniquely as a product $g=nd $, for $n\in \Ni$ and $d\in \Di$. The action by $g$ on $ev\in V$ is determined by
\[g(ev)=\leftexp{d}{e}\cdot\chi(n)\cdot d(v).\]
The endomorphisms of $V$ as a $\ki$-vector space are just $\ki$; the endomorphisms of $V$ that commute with $\Gi$ are exactly those that are fixed by the action of $\Di$, namely, the elements of $\L$.

Now since every object in $\Rep{k}{{\G}}$ has a canonical decomposition as a direct sum of irreducible representations, the result follows.
\end{proof}

\begin{proposition}\label{K0 Repk<G> = Z[Z/linfZ]}
The ring $K_0\RkG$ of isomorphism classes of $k$-semilinear representations of $\G$ is
\[K_0\RkG\cong\ZZ[\QQ_\l/\ZZ_\l].\]
\end{proposition}

\begin{proof}
This follows immediately from Example~\ref{L<Z/pi x Z/pi>} and Proposition~\ref{k<G> is a colimit}.
\end{proof}

Now that we have computed the homotopy groups of $K\Rep{k}{\G}$, we will use them to compute, in Proposition~\ref{htpy grps of dc of Repk<G>}, the homotopy groups of the derived completion $\dc{K\Rep{k}{\G}}{\alpha_\l}$. We will use the following well-known lemma which describes the homotopy groups of the $\l$-adic completion $\dc{K\L}{\l}$.

\begin{lemma}\label{Bousfield completion is tensor Zl}
 Let $X$ be a spectrum with homotopy groups $M_n:=\pi_n X$,
and let $U_n \subseteq M_n$ be a uniquely $\l$-divisible part of $M_n$. Then if $M_n/U_n$ is a finitely generated abelian group, the homotopy groups of the $\l$-adic completion of $X$ are
\[\pi_n(X\hatl) \cong M_n/U_n \otimes \ZZ_\l.\]
\end{lemma}

The next two technical lemmas simplify our computations of the homotopy groups of derived completions.

\begin{lemma}\label{ok to bousfield complete before derived}
Let $\alpha\colon \ZZ[\QQ_\l/\Zl] \to \FF_\l$ be the mod-$\l$ augmentation map and let $X$ be an $H(\ZZ[\QQ_\l/\ZZ_\l])$-module. Let $M_n:=\pi_n (X)$, and define $U_n \subseteq M_n$ to be the uniquely $\l$-divisible part of $M_n$. Then if $M_n/U_n$ is finitely generated, the map $\beta\colon X \to X\hatl$ induces a weak equivalence
\[\big(\ZZ[\QQ_\l/\Zl] \otimes \pi_* X\big)^{^\wedge}_{\alpha} \rTo^{\cong} \dc{\ZZ[\QQ_\l/\Zl] \otimes \pi_* (X\hatl)}{\alpha}\]
on derived completions.
\end{lemma}

\begin{proof}
By Lemma ~\ref{Bousfield completion is tensor Zl},
\[\pi_*(X\hatl) \cong M_n/U_n \otimes \ZZ_\l.\]
There is an isomorphism
\[\FF_\l \otimes M_n \cong \FF_\l \otimes M_n/U_n \otimes \ZZ_\l.\]
and hence, by a theorem of Elmendorf-Kriz-Mandell-May~\cite{EKMM}, $\beta$ induces weak equivalences between the $n$-simplices of the derived completions.
\end{proof}

\begin{lemma}\label{finitely generated abelian pulls out of dc}
Let $R \to S$ be a map of commutative rings and $M$ an $R$-module with no $\ZZ$-torsion. Then if $A$ is a finitely generated abelian group,
\[\pi_* \left(\dc{M \otimes A}{S}\right) \cong \pi_*(M^{^\wedge}_S) \otimes A.\]
\end{lemma}

\begin{proof}
Let $0 \to \ZZ^n \rTo^{i} \ZZ^m \to A \to 0$ be a free resolution of $A$. The result now follows immediately from the fact that the derived completion and $\pi_*$ preserve fiber sequences.
\end{proof}

\begin{proposition}\label{htpy grps of dc of Repk<G>}
The algebraic to geometric spectral sequence for $\dc{K\RkG}{\alpha_\l}$ collapses at the $E_1$-term and thus $\dc{K\RkG}{\alpha_\l}$ has homotopy groups
\[\pi_i\left(\dc{K\RkG}{\alpha_\l}\right) \cong \pi_i \left((K\FF_{p})\hatl\right) \oplus \pi_{i-1}\left((K\FF_{p})\hatl\right).\]
\end{proposition}

\begin{proof}
The $E_1$-term of the algebraic to geometric spectral sequence is
\[E_1^{s,t}=\pi_{s+2t}\left(\dc{\ZZ[\QQ_\l/\Zl] \otimes \pi_{-t} K\L\phantom{\widehat{A}}\hspace{-.13in}}{\FF_\l}\right).\]
By Lemmas \ref{Bousfield completion is tensor Zl} and \ref{finitely generated abelian pulls out of dc} along with Example~\ref{dc of Z/linfZ}, this becomes
\[
E_{1}^{s,t} = \left\lbrace
\begin{array}{ll}
\pi_{-t} \left((K\L)\hatl\right)\quad & \text{when }s+2t=0\text{ or }1 \\
0 & \text{otherwise}.
\end{array}
\right.
\]
Our $E_1$ page is therefore as in the diagram below.\\



        \[ \xygraph{
!{<0cm,0cm>;<2.22cm,0cm>:<0cm,1.5cm>::} 
    !{(0,0)}*+{t}="00" 
    !{(1,0)}*{0}="10"
    !{(2,0)}*{0}="20"
    !{(3,0)}*{0}="30"
    !{(4,0)}*{0}="40"
    !{(5,0)}*{\pi_0(\dc{K\L}{\l})}="50"
    !{(0,1)}="01"
    !{(1,1)}*{0}="11"
    !{(2,1)}*{0}="21"
    !{(3,1)}*{0}="31"
    !{(4,1)}*{0}="41"
    !{(5,1)}*+{\pi_0(\dc{K\L}{\l})}="51"
    !{(0,2)}="02"
    !{(1,2)}*{0}="12"
    !{(2,2)}*{0}="22"
    !{(3,2)}*{0}="32"
    !{(4,2)}*+{\pi_1(\dc{K\L}{\l})}="42"
    !{(5,2)}*{0}="52"
    !{(0,3)}="03"
    !{(1,3)}*{0}="13"
    !{(2,3)}*{0}="23"
    !{(3,3)}*{0}="33"
    !{(4,3)}*+{\pi_1(\dc{K\L}{\l})}="43"
    !{(5,3)}*{0}="53"
    !{(0,4)}="04"
    !{(1,4)}*{0}="14"
    !{(2,4)}*{0}="24"
    !{(3,4)}*+{\pi_2(\dc{K\L}{\l})}="34"
    !{(4,4)}*{0}="44"
    !{(5,4)}*{0}="54"
    !{(0,5)}="05"
    !{(1,5)}*{0}="15"
    !{(2,5)}*{0}="25"
    !{(3,5)}*+{\pi_2(\dc{K\L}{\l})}="35"
    !{(4,5)}*{0}="45"
    !{(5,5)}*{0}="55"
    !{(5,6)}*+{s}="56"
    !{(2.3,5.4)}="255" 
    !{(2.8,5.4)}="2854" 
    "50":"00" 
    "50"-"51" "51":"56" 
    "50"-@{.}"42" "42"-@{.}"34" "34"-@{.}"255" 
    "51"-@{.}"43" "43"-@{.}"35" "35"-@{.}"2854" 
    }\]



\noindent By Proposition \ref{l-comp(KL)=l-comp(KFl)}, the
diagonal lines each correspond to a copy of
$\pi_*\big((K\FF_p)\hatl\big)$.  All of the differentials vanish, and hence the spectral sequence collapses. Therefore the homotopy groups of the derived completion are
\[ \pi_i\left(\dcl{K\RkG}\right) \cong \pi_i\left((K\FF_{p})\hatl\right) \oplus \pi_{i-1}\left((K\FF_{p})\hatl\right).\]
\end{proof}

\section{The map $\E$}
\label{weak equiv}

In the proof of Theorem~\ref{maintheorem}, we will follow a strategy similar to the one found in Section~\ref{homotopy groups} of \cite{CarlDRT}. Define $k_i$, $\G_i'$ and $\Gi$ as in example~\ref{L<Z/pi x Z/pi>}, let $\F_i:=(\F^t)^{\G_i'}$ and let $\O[i]$ be the ring of integers in $\Fi$. Denote by $\ROi$ the category of finitely generated $\O[i]$-modules equipped with the continuous semilinear action by $\Gi$ that is inherited from the natural action by $\Gi$ on $\Fi$.

The categories $\colim \ROi$ and $\colim \Rep{\Fi}{\Gi}$ are
closed under tensor products with objects in
$\RkG$, and thus their $K$-theory
spectra are both $K\RkG$-modules. We may therefore take their
derived completions across the augmentation map $\alpha_\l:
K \RkG \to K \mathfrak{F}_\ell$
(equation~\ref{aug}). As in \cite{CarlDRT}, we will show that the
intermediate category $\colim \ROi$ satisfies the two conditions
below. These conditions will be proven in Propositions \ref{colim
  KRep_O<G_i> = KRep_k<G>} and \ref{KRep_O<G_i> = KRep_Fi<G_i>},
respectively.
\begin{enumerate}
\item Extension of scalars induces a weak equivalence of completed $K$-theory spectra,
\[\bigdca{K\RkG}\xto{\simeq}\bigdca{\colim K\ROi}.\]

\item Extension of scalars induces a weak equivalence of completed $K$-theory spectra,
\[\bigdca{\colim K\ROi}\xto{\simeq} \bigdca{K\Rep{\F^t}{\G}}\simeq \dcl{K\F}.\]
\end{enumerate}

We start by proving a fact about $K$-theory spectra of rings of
the form $R =E [[x^\frac1r]] \langle G \rangle$ where $E$ is a
field, $G$ is a finite group, and $r$ is a positive
integer.  Denote by  $\Mod {R}$ the category of finitely generated
projective $R$-modules, and by  $\NoTor{R} \subseteq \Mod {R}$ the
full subcategory of finitely generated torsion free projective
$R$-modules.

\begin{lemma}\label{Tors free K-theory is w. eq.}
The inclusion functor $\NoTor{R} \hookrightarrow \Mod {R}$ induces a weak equivalence of $K$-theory spectra.
\end{lemma}

\begin{proof}
Let $M$ be a projective $R$-module on $n$ generators, where $n$ is a positive integer. Then we can construct a surjective map $R^n \to M$. The kernel of this map is torsion free, and thus we have a length 2 resolution of $M$ by torsion free modules. The result now follows by ``Reduction by Resolution," Theorem 3.3 of \cite{Q}.
\end{proof}

\begin{proposition}\label{colim KRep_O<G_i> = KRep_k<G>}
Extension of scalars induces a weak equivalence,
\[\bigdca{\colim K\RkGi}\xto{\simeq} \bigdca{\colim K\ROi}.\]
\end{proposition}

\begin{proof}
First we show that the map $t_i:\RkGi\xto{\Oi\bracket{\Gi}\tensorover{\ki\bracket{\Gi}}(-)}\ROi$ induces an isomorphism on $\pi_0$.

The categories $\RkGi$ and $\ROi$ can be identified with the categories $\Mod {\ki \bracket{\Gi}}$ and $\Mod {\Oi \bracket{\Gi}}$, respectively, since the group $\Gi$ is finite.  By Lemma \ref{Tors free K-theory is w. eq.}, it suffices to show that $t_i$ induces an isomorphism
\[\tilde{t}_i\colon\pi_0 K\NoTor{\ki \bracket {\Gi}} \xto{\simeq} \pi_0 K\NoTor{\Oi \bracket {\Gi}}.\]

Clearly, $\tilde{t}_i$ is an injection on the set of isomorphism classes of objects. To see that it is a surjection, we start with a module $M$ over the ring $\Oi \bracket{\Gi}$. Let $I$ be the kernel of the quotient map $\Oi \bracket{\Gi} \to \ki\bracket{\Gi}$; $I$ is generated by the fractional powers of $x$. The module $IM$ is closed under multiplication by $\Oi \bracket{\Gi}$, and hence the quotient map $M \xto{\gamma} M/IM$ is a map of $\Oi \bracket{\Gi}$-modules. Both $M$ and $M/IM$ are, by restriction, modules over $\ki \bracket{\Gi}$ as well, and since $\ki \bracket{\Gi}$ is semisimple, the map $\gamma$ splits. Therefore we have a section $s\colon M/IM \to M$. This map can be extended to a map of $\Oi \bracket{\Gi}$-modules,
\[ S\colon\Oi \bracket{\Gi}  \tensorover{\ki \bracket{\Gi}} M/IM \to M.\]

Now as a consequence of Nakayama's Lemma, a homomorphism $f\colon X \to Y$ of free modules over a complete local ring with maximal ideal $M$ is an isomorphism exactly when the induced map on the quotient is.

Since $M$ and $\Oi \bracket{\Gi}  \tensorover{\ki \bracket{\Gi}} M/IM$ are free over $\Oi$ and the induced map on quotients
\[\big(\Oi \bracket{\Gi}  \tensorover{\ki \bracket{\Gi}} M/IM\big) \Big/ I\big(\Oi \bracket{\Gi}  \tensorover{\ki \bracket{\Gi}} M/IM\big) \to M/IM\]
is an isomorphism, the map $S$ must be an isomorphism of $\Oi$-modules. Thus every module over $\Oi \bracket{\Gi}$ is ``extended" from $\ki \bracket{\Gi}$.

Next we show that the endomorphism rings of the irreducible modules over $\Oi\bracket{\Gi}$ are of the form $\L[[x^{(1/p)},x^{(1/\l^{i+1})}]]$ by examining the irreducible $\ki\bracket{\Gi}$-modules.

By Example~\ref{L<Z/pi x Z/pi>} the irreducible $k_i [\Gi]$-modules are in bijective correspondence with the group $\Hom(\N_i, k)$, where $\N_i=\Zmod{\l^i}$. Suppose that $M_\alpha$ is an irreducible $\ki \bracket{\Gi}$ module corresponding to a character $\alpha \in \Hom(\Ni,k)$. Then because $\Gi$ acts on $M_{\alpha}$ through $\alpha$, the tensor product $M_{\alpha}':=\Oi\bracket{\Gi}\tensorover{\ki\bracket{\Gi}}M_{\alpha}$ is isomorphic as an $\O\langle\Gi\rangle$-module to $\Oi$. The endomorphisms of $M_\alpha'$ over $\Oi \bracket{\Gi}$ are exactly the elements of $M_\alpha'$ that commute with the action by $\Oi\langle\Gi\rangle$---namely, $\L[[x^{(1/p)},x^{(1/\l^{i+1})}]]$.

The homotopy groups of $K\ROi$ are therefore
\[ \pi_{*}K\ROi = K_0\ROi \otimes K_* \L[[x^{(1/p)},x^{(1/\l^{i+1})}]],
\]and thus the homotopy groups of the colimit are
\[\pi_*\left(\colim K\ROi\right) = \ZZ[\QQ_\l/\Zl] \otimes K_* \L[[x^{(1/p)}]].\]

We now have the commutative diagram below. We set $\mathfrak{Z}:=\ZZ[\QQ_\l/\Zl]$ in all of the diagrams of this proof for typesetting purposes.
\begin{diagram}
\pi_*K\RkG\cong\mathfrak{Z} \otimes \pi_*K\L&\rTo&\mathfrak{Z} \otimes \pi_*\big((K\L)\hatl\big)\\
\dTo_{t_*}& &\dTo_{(t_*)\hatl}\\
\pi_*\colim K\ROi \cong\mathfrak{Z}\otimes \pi_* K\L[[x^{(1/p)}]] & \rTo & \mathfrak{Z} \otimes \pi_*\big((K\L[[x^{(1/p)}]])\hatl\big)\\
\dTo_{q}& &\dTo_{(q)\hatl}^{\cong}\\
\pi_*K\RkG \cong \mathfrak{Z}\otimes K\L & \rTo & \mathfrak{Z} \otimes \pi_*\big((K\L)\hatl\big)\\
\end{diagram}
\noindent The horizontal arrows are induced by the natural maps that come from the Bousfield completion, the map $t$ is defined to be the colimit of the $t_i$, and the map $t_*$ is the induced map on homotopy groups. The map $q$ is induced by the quotient map $\L[[x^{(1/p)}]] \to \L$ and, by a theorem of Suslin~\cite{Sus84}, induces a weak equivalence on the $\l$-adic completions. The composition $q \circ t$ is the identity and hence if $\dc{q}{\l}$ is a weak equivalence, $(t_*)\hatl$ must also be a weak equivalence. The maps in the top two rows of the diagram above induce the diagram of $E_1$-terms for the algebraic to geometric spectral sequence below.
\begin{diagram}
\pi_{s+2t}\left(\dc{\mathfrak{Z} \otimes \pi _{-t}K\L}{\FF_\l}\right)
& \rTo^{\cong}
& \pi_{s+2t}\left(\sdc{\left(\mathfrak{Z} \otimes \pi_{-t}K\L\hatl\right)}{\FF_\l}\right)\\
\dTo_{\varepsilon(t_*)}& & \dTo_{\varepsilon(t_*)\hatl}\\
\pi_{s+2t}\left(\big(\mathfrak{Z} \otimes \pi_{-t} \sdc{K\L[[x^{(1/p)}]]\big)}{\FF_\l}\right)
& \rTo^{\cong}
& \pi_{s+2t}\left(\left(\mathfrak{Z} \otimes \pi_t \sdc{K\L[[x^{(1/p)}]]\hatl \right)}{\FF_\l}\right).
\end{diagram}

By Lemma~\ref{ok to bousfield complete before derived} the horizontal maps are isomorphisms on $E_1$-terms. Therefore the map $\varepsilon(t_*)$ of $E_1$-terms is an isomorphism. Since by Lemma ~\ref{htpy grps of dc of Repk<G>}, the spectral sequence collapses at $E_1$, the map $t$ induces a weak equivalence,
\[\bigdca{K\RkG} \xto{\simeq} \bigdca{\colim K\ROi}.\]
\end{proof}

\begin{proposition}\label{KRep_O<G_i> = KRep_Fi<G_i>}
Extension of scalars induces a weak equivalence of completed $K$-theory spectra, \[\dcl{\colim K\ROi}\xto{\simeq}\dcl{K\Rep{\F^t}{\G}} \cong \sdcl{K\F}.\]
\end{proposition}
\begin{proof}
As previously mentioned, the weak equivalence
$\dcl{K\Rep{\F^t}{\G}} \simeq \dcl{K\F}$ is an immediate
consequence of the equivalence of categories $\Rep{\F^t}{\G}\cong
\Mod{\F}$.  Therefore we need only show that \[\dcl{\colim K\ROi}\to \dcl{K\Rep{\F^t}{\G}}\] is a weak equivalence. Since $\Gi$ is finite, the category $\ROi$ can also be thought of as a category of finitely generated modules over the skew group ring $\O[i]\langle \Gi\rangle$. In this ring, inverting all powers of $x$ (including fractional ones) yields the ring $\Fi$, and hence we can use Corollary 5.12 from \cite{Swan} and Quillen's localization theorem \cite{Q} to conclude that there is a homotopy fiber sequence of $K$-theory spectra,
\[K\Ti\to K\ROi\to K\Rep{\Fi}{\Gi},\]
where $\Ti$ is the full subcategory of modules in $\ROi$ that localize to $0$. These are precisely the modules that are annihilated by some (possibly fractional) power of $x$.

The category $\colim \Ti$ is closed under the tensor product with
objects in $\RkG$. We may therefore take its derived completion
across the map $\alpha_\l$ (see equation~\ref{aug}). Since
colimits and derived completions preserve
fiber sequences, we have a fiber sequence of $K\RkG$-modules
\[\dcl{\colim K\Ti} \to \dcl{\colim K\ROi} \to \dcl{\colim  K\Rep{\Fi}{\Gi}}.\]
Therefore it suffices to show that $\dcl{\colim K\Ti}$ is contractible.

Let $\T:=\colim \Ti$, let $\Ri:=\RkGi$, and let $\R:=\colim\Ri$. The derived completion $\KTl$ is the cosimplicial $\KR$-module whose $n$-simplices are obtained by $n+1$ applications of the functor $(\Hl\smashover{\KR}-)$ to $\KT$. It suffices to show that $\Hl\smashover{\KR}\KT$, which is isomorphic to $\colim (\Hl\smashover{\KRi}\KT_i)$, is contractible.\label{KT smash Hl is a colim}

We will now show that for every $i$, the $K$-theory spectrum $\KTi$ is weakly equivalent to $\KRi$. Let $\XX \subset \NN$ be the subcategory of positive integers whose reciprocal occurs as an exponent of $x$ in the ring $\Oi$. For every $s \in \XX$ define $\Zi[s]\subset\Ti$ to be the full subcategory of modules that are killed by $x^{1/s}$. Then by devissage~\cite{Q}, the $K$-theory spectrum $K\Zi[s]$ is weak equivalent to $K\Ti$. On the other hand, the category $\Zi[s]$ is equivalent to the category $\Rep{\Oi{[x^{(1/s)}]}}{\Gi}$ via the forgetful functor.

Because for every positive integer $r$, the category $\Zi[rs]$ is a full subcategory of $\Zi[s]$, the functor $\Zi[-]$ takes the category $\XX$ contravariantly to the category of full subcategories of $\Ti$; it is therefore a cofiltered diagram. The category $\Ri$ maps via $\tensorup{\Oi{[x^{(1/s)}]}}{\ki}$ to $\Rep{\Oi{[x^{(1/s)}]}}{\Gi}$, and it is easy to see that it is the limit
\[\Ri=\lim_{s\in\XX}\Rep{\Oi{[x^{(1/sr)}]}}{\Gi}.\]
Applying the functor $K$ yields the commutative diagram
\begin{diagram}
K\Z_i&\rTo&\ldots&\rTo^{\simeq}&K\Zi[sr]&\rTo^{\simeq}&K\Zi[s]&\rTo^{\simeq}K\Ti\\
\dTo_{\simeq}&    &      &            &  \dTo_{\simeq} &            & \dTo_{\simeq} &            &\\
\KRi&\rTo&\ldots&\rTo&K\Rep{\Oi{[x^{(1/sr)}]}}{\Gi}&\rTo&K\Rep{\Oi{[x^{(1/s)}]}}{\Gi}
\end{diagram}
where $\Z_i$ is defined to be the limit $\lim_{s\in \NN}\Zi[s]$. This gives us the weak homotopy equivalence,
\[K\Ti\simeq K\Z_i\simeq K\Ri.\]

Now, since $\pi_*$ commutes with filtered colimits and smash products,

the above equivalence yields
\begin{align}
\pi_*(\Hl\smashover{\KR}\KT)\cong\colim\pi_*(\Hl\smashover{\KRi}\KZi)\cong\colim\pi_*\Hl.\label{eqtn-pi* KT smash Hl = colim pi*Hl}
\end{align}
Since $\pi_n\Hl=0$ for all $n\neq 0$, we restrict our attention to $\pi_0$.
Because $\Z_i$ is closed under tensor product, $K\Z_i$ is actually a ring spectrum, so $\pi_0K\Z_i$ is a ring, and
\[\pi_0(\Hl\smashover{K\Ri}K\Z_i)=\pi_0\Hl\tensorover{}\pi_0K\Z_i,\] where the tensor product is taken over $\pi_0K\Ri$.
As a free $\pi_0K\Ri$-module, $\pi_0K\Z_i$ is generated by a
$1$-dimensional $\ki$-vector space $V_i$ that is annihilated by
all (possibly fractional) powers of $X$ in $\Oi$. Under
$\tau_i^j$, the module $V_i$ (which corresponds to $1$ in
$\pi_0\Hl$) gets sent to the module
\[V_i \tensorover{\Oi} \O[j]=\ki[j]\big[\big[x^{(1/\l)},x^{(1/p)}\big]\big]\,[x^{1/{\l^j}}]\]
(which corresponds to $\l^j$ in the ring
$\pi_0\Hl=\FF_\l$), and thus
all of the maps in the colimit equation~(\ref{eqtn-pi* KT smash
  Hl = colim pi*Hl}) are zero maps. Therefore $\pi_*(\Hl
\smashover{K\R} K\T)$ is contractible.

The $n$-simplices of $\dcl{\KT}$ are of the form $\Hl \smashover{K\R} ... \smashover{K\R} \Hl \smashover{K\R} K\T$ so they are all contractible as well. Since a level-wise weak equivalence of cosimplicial objects induces a weak equivalence on total spaces, $\sdcl {K\T}$ is contractible, and thus
\[\dcl{\colim K\ROi}\xto{\simeq} \dcl{K\Rep{\F^t}{\G}}.\]\end{proof}
The following theorem follows immediately from the two preceding propositions.

\begin{theorem}\label{maintheorem}
Extension of scalars induces a weak equivalence of completed $K$-theory spectra,
\[\dcl{K\RkG}\xto{\simeq}\dcl{K\Rep{\F^t}{\G}}\cong\dcl{K\F}.\]
\end{theorem}

\newpage
\bibliographystyle{alpha}
\bibliography{bibliography}
\begin{table}[b]
\center Department of Mathematics, Massachusetts Institute of Technology\\
Cambridge, MA 02139\\
gracelyo@math.mit.edu
\end{table}
\end{document}